\documentclass[12pt]{article}
 \usepackage{latexsym,amsmath,amssymb,amsfonts}
 \usepackage[mathscr]{eucal}
 \usepackage[all]{xypic}
\author{Charles Cadman}
\title{{\bf On the enumeration of rational plane curves with tangency conditions}}
\date{}

\newtheorem{contact_type}{Proposition}[section]
\newtheorem{nice_maps}[contact_type]{Definition}
\newtheorem{inf_def}[contact_type]{Lemma}
\newtheorem{Cap_Harris_1}[contact_type]{Lemma}
\newtheorem{Cap_Harris_2}[contact_type]{Lemma}
\newtheorem{gen_smooth}[contact_type]{Theorem}

\newtheorem{relations}{Proposition}[section]
\newtheorem{CH_numbers}[relations]{Definition}
\newtheorem{CH_enum}[relations]{Proposition}
\newtheorem{CH_recursion}[relations]{Theorem}

\newtheorem{key_relation}{Proposition}[section]

\newtheorem{enumerative}{Theorem}[section]
\newtheorem{smooth_source}[enumerative]{Lemma}
\newtheorem{gen_elmt}[enumerative]{Lemma}

\newcommand{\Pic}{\mathrm{Pic}}

\newcommand{\Der}{\mathrm{Der}}

\newcommand{\spec}{\mathrm{Spec}\;}

\newcommand{\mult}{\mathrm{mult}}

\newcommand{\edim}{\mathrm{edim}}

\newcommand{\tors}{\mathrm{tors}}

\newcommand{\sO}{\mathcal{O}}

\newcommand{\sL}{\mathcal{L}}

\newcommand{\sN}{\mathcal{N}}

\newcommand{\sU}{\mathcal{U}}

\newcommand{\sX}{\mathcal{X}}

\newcommand{\fC}{\mathfrak{C}}

\newcommand{\ff}{\mathfrak{f}}

\newcommand{\eC}{\mathscr{C}}
\newcommand{\eK}{\mathscr{K}}
\newcommand{\eU}{\mathscr{U}}

\newcommand{\bmu}{\boldsymbol{\mu}}
\newcommand{\bZ}{\mathbb{Z}}

\newcommand{\bC}{\mathbb{C}}

\newcommand{\bP}{\mathbb{P}}

\newcommand{\bI}{\mathbb{I}}

\newcommand{\pd}{\bP^2_{D,r}}
\newcommand{\pe}{\bP^2_{E,2}}
\newcommand{\beginproof}{\emph{Proof:} }

\newcommand{\qed}{\hfill $\Box$\smallskip}

\newcommand{\ch}[2]{\genfrac{(}{)}{0pt}{}{#1}{#2}}

\begin{document}
\maketitle
\begin{abstract}
We use twisted stable maps to answer the following question.  Let $E\subset\bP^2$ be a smooth cubic.  How many rational degree $d$ curves pass through $a$ general points of $E$, have $b$ specified tangencies with $E$ and $c$ unspecified tangencies, and pass through $3d-1-a-2b-c$ general points of $\bP^2$?  The answer is given as a generalization of Kontsevich's recursion.  We also investigate more general enumerative problems of this sort, and prove an analogue of a formula of Caporaso and Harris.
\end{abstract}

\section{Introduction}
\label{sec:intro}
\setcounter{equation}{0}

It is well-known that Gromov-Witten theory can be used to answer questions in enumerative geometry.  This dates back to Kontsevich's celebrated recursion for the number $N_d$ of rational plane curves passing through $3d-1$ general points (here $k,\ell\ge 1$).
$$N_d=\sum_{k+\ell=d} N_kN_{\ell} k^2\ell\left[\ell\ch{3d-4}{3k-2}-k\ch{3d-4}{3k-1}\right],\; d\ge 2$$
In this paper, we follow up on the results of \cite{stacks,GWinvs} by showing that the Gromov-Witten invariants of a certain Deligne-Mumford stack, $\bP^2_{E,2}$, are enumerative.  Here $E\subset\bP^2$ is a smooth cubic, and $\bP^2_{E,2}$ is obtained from $\bP^2$ by applying the square root construction along $E$, which was defined in \cite{stacks}.  The main property of this stack which makes it interesting for enumerative geometry is that a morphism $f:C\to\bP^2$ from a smooth connected curve $C$ such that $f(C)\not\subset E$ lifts to a morphism $C\to\bP^2_{E,2}$ if and only if $f^*E$ has even multiplicity at every point.  This leads to the idea that the stack of twisted stable maps to $\bP^2_{E,2}$ looks like a stack of maps to $\bP^2$ with tangency conditions imposed.  This was verified in \cite{stacks}.

In \cite{CH}, Caporaso and Harris found a recursion which computes the number of plane curves of genus $g$ passing through a specified number of points and having certain contact conditions with respect to a line.  These contact conditions come in the form of $k$th order contacts at specified points of the line or at unspecified points.  This was generalized by Vakil \cite{Va}, who solved the corresponding problem with the line replaced by a smooth rational curve on a rational surface (see [ibid] for hypotheses).  In particular, he solved the problem for contacts with a smooth plane conic.  This paper is a step towards solving the problem for a smooth plane cubic $E$.

We answer the following question.  Given integers $a,b,c,d$ satisfying appropriate conditions for the question to be sensible, how many rational degree $d$ curves meet $E$ in $a$ specified order $1$ contacts, $b$ specified order $2$ contacts, and $c$ unspecified order $2$ contacts, and pass through $3d-1-a-2b-c$ general points in $\bP^2$.  Here by a specified contact, we mean that we have chosen a general point of $E$ and require the curve to pass through this point and have a certain order of contact there.  By an unspecified contact, we mean that the contact occurs at an arbitrary point distinct from the specified points.  We require all these contacts to occur at smooth points of the rational curve.  We show that for general choices of points, the rational curves satisfying these conditions have at worst nodal singularities, none of which lie on $E$.  In section \ref{sec:kontsevich} we denote this number by $N_d(a,b,c)$, and the solution is given in equations \ref{eq:main_1}-\ref{eq:main_3}.  A maple program implementing this recursion can be downloaded from the author's homepage.

In section \ref{sec:def}, we begin by recalling facts that we proved elsewhere.  Then we use some results of Caporaso and Harris to prove some nice facts about the stack of twisted stable maps into $\bP^2_{D,r}$.  Finally, we discuss the relation between deformations of a map from a twisted curve having a separating node with deformations of the maps from the two curves obtained by normalizing the node.  This is essential to prove enumerativity of the Gromov-Witten invariants.

In section \ref{sec:cap_harris}, we show that the number of curves with specified contacts to a smooth curve in $\bP^2$ equals an intersection number on the stack of twisted stable maps.  Then we give a relation between these numbers which was inspired by the recursion of Caporaso and Harris \cite{CH}.  In fact, it is the formula one would expect to get if the curves being deformed were not allowed to break into pieces.  In our case they cannot, because a rational curve cannot map onto a curve of positive genus.

In section \ref{sec:kontsevich}, we derive the main result modulo the proof of enumerativity, which is given in section \ref{sec:enumerative}.

Several results in this paper can easily be generalized, but in the interest of efficiency we restrict attention to plane curves.

\smallskip

\noindent{\bf Notation and Conventions}

As is common in enumerative geometry, we work over $\bC$.  Throughout this paper, $D$ denotes a smooth plane curve of degree $\delta\ge 3$.
\smallskip

\noindent{\bf Acknowledgements}

The fact that certain twisted Gromov-Witten invariants are enumerative is part of my Ph.D. thesis at Columbia University under the supervision of Michael Thaddeus.  This paper extends the enumerative result from my thesis by allowing tangencies at arbitrary points, and also verifies a conjecture contained therein.

\section{Deforming morphisms from twisted curves into $\bP^2_{D,r}$}
\label{sec:def}
\setcounter{equation}{0}

In this section we review some facts about twisted stable maps into $X_{D,r}$ and then study their deformations.  When the domain is smooth and does not map into $D$, a twisted stable map is equivalent to an ordinary map with tangency conditions to $D$ imposed.  Much of what we need to understand deformations in this case is already worked out in \cite{CH}, and we review their results.  We also study deformations of a twisted stable map preserving a node.

\subsection{Review of twisted stable maps}

Twisted stable maps were defined by Abramovich and Vistoli in \cite{AV}.  We are interested in twisted stable maps into a particular kind of stack.  Let $D\subset\bP^2$ be a smooth curve of degree $\delta\ge 3$ and let $r$ be a positive integer.  In \cite{stacks}, we introduced the $r$th root construction, which produces a Deligne-Mumford stack $\pd$.  Locally it is the quotient of a cyclic $r$ sheeted covering of $\bP^2$ which it totally ramified along $D$ by the $\mu_r$-action \cite[2.15]{stacks}.  In \cite[\S 2.1]{GWinvs}, it is shown that there is a discrete invariant on the stack of twisted stable maps into $\pd$ called the contact type.  If $n$ is the number of marked points, then the contact type is an $n$-tuple of integers between $0$ and $r-1$.  The contact type allows for a nice characterization of maps from smooth twisted curves into $\pd$.

The following proposition is a consequence of \cite[3.9]{stacks}.  Note that since we assumed $\delta\ge 3$, no rational curve can map onto $D$.

\begin{contact_type}
\label{th:contact_type}
Let $\fC$ be a smooth, $n$-marked, genus $0$ twisted curve over a scheme $S$ and let $\ff:\fC\to\pd$ be a twisted stable map of positive degree and contact type $\vec{\varrho}$.  Let $C$ be the coarse moduli space of $\fC$ with induced markings $\sigma_i\subset C$, and let $f:C\to\bP^2$ be induced by $\ff$.  Then there is an effective Cartier divisor $Z\subset C$ such that
\begin{equation}
\label{eq:contact_cond}
f^*D=rZ+\sum_{i=1}^n\varrho_i\sigma_i.
\end{equation}
Moreover, given a morphism $f:C\to\bP^2$ and an effective Cartier divisor $Z\subset C$, there is a unique (up to isomorphism) twisted curve $\fC$ with coarse moduli space $C$ and a unique twisted stable map $\ff:\fC\to\pd$ with contact type $\vec{\varrho}$ which induces $f$.
\end{contact_type}

The expected dimension of the stack $\eK_{0,n}(\pd,d,\vec{\varrho})$ is
\begin{equation}
\label{eq:edim}
d(3-\delta)+\frac{1}{r}(d\delta - \sum \varrho_i) + n - 1.
\end{equation}
This was computed in equation 3.7 of \cite{GWinvs}.

\begin{nice_maps}
\label{th:nice_maps}
We define $\eK^{*}_{0,n}(\pd,d,\vec{\varrho})$ as follows.  Define an irreducible component of $\eK_{0,n}(\pd,d,\vec{\varrho})$ to be \emph{good} if the general map in that component has a smooth source curve and maps birationally onto its image.  Let $\eU\subset\eK^{*}_{0,n}(\pd,d,\vec{\varrho})$ be the open substack obtained by removing all the bad irreducible components, and let $\eK^*_{0,n}(\pd,d,\vec{\varrho})$ be its stack-theoretic closure.
\end{nice_maps}

Note that $\eK^*_{0,n}(\pd,d,\vec{\varrho})$ has an open dense substack which is a scheme.  Indeed, the general map in each irreducible component has no automorphisms, and $\eK_{0,n}(\pd,d)$ has a projective coarse moduli scheme \cite[1.4.1]{AV}.

\subsection{Maps from smooth twisted curves}

We begin with an easy lemma about certain types of infinitesimal deformations.  Let $X$ be a nonsingular curve, let $Y$ be a nonsingular surface, let $f:X\to Y$ be a morphism whose differential $df:T_X\to f^*T_Y$ is an injection of sheaves, and let $D\subset Y$ be a nonsingular curve.  Let $\sN$ be the cokernel of $df$, and let $f^*D=\sum_{i=1}^n \varrho_ip_i$, where $p_i\in C$ are distinct points and $\varrho_i>0$.  Assume that $df$ is injective on fibers at each $p_i$.

\begin{inf_def}
\label{th:inf_def}
The first order infinitesimal deformations of $f$ which fix $Y$ (but not necessarily $X$) and preserve the multiplicities $\varrho_i$ (but not the points $x_i$) are naturally in bijection with $H^0(X,\sN(-\sum (\varrho_i-1)x_i))$.
\end{inf_def}

\beginproof  First we recall the construction which identifies first order deformations of $f$ fixing $Y$ with $H^0(X,\sN)$.  Let $\bI=\spec k[\epsilon]/(\epsilon^2)$, and suppose we have the following commutative diagram.
$$\xymatrix{
 & Y \\
X \ar[ur]^f \ar@{^{(}->}[r] \ar[d] \ar@{}[dr]|{\Box} & \sX \ar[d] \ar[u]_F \\
\spec k \ar@{^{(}->}[r] & \bI}$$
Cover $X$ by affines $U_i$.  Since $X$ and $\sX$ have the same underlying topological space, we get an open covering of $\sX$ by affines $\sU_i$ so that $U_i$ embeds into $\sU_i$ as $(\sU_i)_{\mathrm{red}}$.  Since nonsingular affine varieties have no nontrivial first order deformations, there are isomorphisms $\varphi_i:U_i\times\bI\to\sU_i$.  On the overlaps $U_{ij}$, $\varphi_j^{-1}\circ\varphi_i$ determines a derivation $\alpha_{ij}\in H^0(U_{ij},T_{U_{ij}})$.  These form a 1-cycle, and hence determine an element $\alpha\in H^1(X,T_{X})$.

On $U_i$, the morphisms $F$ and $\varphi_i$ determine a morphism $\psi_i:U_i\times\bI\to Y$, and hence a derivation $\beta_i\in H^0(X,f^*T_Y\vert_{U_i})$.  The commutative diagram
$$\xymatrix{
U_{ij}\times\bI \ar[r]^{\varphi_i} \ar[dr]_{\psi_i} & \sU_{ij} \ar[r]^{\varphi_j^{-1}} & U_{ij}\times\bI \ar[dl]^{\psi_j} \\
 & Y &}$$
shows that on $U_{ij}$ we have $\beta_i=df(\alpha_{ij})+\beta_j$.  It follows that the $\beta_i$ glue to give an element $\beta\in H^0(X,\sN)$.  This is independent of the choices and identifies the first order deformations with $H^0(X,\sN)$.

We need a necessary and sufficient condition for the first order deformation corresponding to $\beta$ to preserve the multiplicities of $f^*D$.  For this we can reduce to the affine situation $U_i=\spec S\to \spec R\subset Y$ and assume that $f^*D=\varrho p$ for some $p\in U_i$.  Let $r\in R$ be a local equation for $D$ and $s\in S$ a local equation for $p$ so that $f^*r=us^{\varrho}$ with $u\in S$ a unit.  The derivation $\beta_i\in\Der_{\bC}(R,S)$ sends $r$ to $f^*r+\beta_i(r)\epsilon$.  This defines a multiplicity $\varrho$ divisor if and only if there are elements $v,w,x,y\in S$ with $v$ a unit such that $$us^{\varrho}+\beta_i(r)\epsilon = (v+w\epsilon)(x+y\epsilon)^{\varrho}.$$  This is equivalent to $\beta_i(r)$ being an element of the ideal generated by $s^{\varrho-1}$.  Now assume $\varrho>1$ since otherwise no conditions are being imposed.  Since $X$ is a curve, $Y$ is a surface, and $df$ is injective at $p$, the condition on $\beta_i$ says precisely that the image of $\beta_i$ in $H^0(U_i,\sN\vert_{U_i})$ vanishes to order at least $\varrho-1$ at $p$.  This finishes the proof.\qed

Proposition \ref{th:contact_type} implies that the deformation theory of twisted stable maps $\fC\to\pd$ from smooth twisted genus $0$ curves $\fC$ is equivalent to the deformation theory of maps $C\to\bP^2$ from smooth rational marked curves $C$ with contact conditions imposed at the markings.  This explains the importance of the above lemma, as well as the following two results of Caporaso and Harris.

Let $\pi:C\to B$ be a smooth, proper family of connected curves over a smooth base $B$, let $f:C\to\bP^2$ be a morphism, and let $b\in B$ be a general point.  Assume that no fiber of $\pi$ maps to a point under $f$.  Let $\sN_b$ be the cokernel of the differential $df_b:T_{C_b}\to f_b^*T_{\bP^2}$, which is injective by hypothesis.  We have a morphism $\kappa_b:T_b B\to H^0(C_b,\sN_b)$ induced by the family of morphisms $C\to B\times\bP^2$, and this is often called the Horikawa map in light of Horikawa's foundational work \cite{Ho}.

\begin{Cap_Harris_1}
\label{th:Cap_Harris_1}
{\bf \cite[2.3]{CH}}  Let $b\in B$ be a general point and assume that $f_b$ maps $C_b$ birationally onto its image.  Then $$\mathrm{Im}(\kappa_b)\cap H^0(C_b,(\sN_b)_\tors)=0.$$
\end{Cap_Harris_1}

Let $Q\subset C$ be the image of a section of $\pi$ such that $f^*D$ has multiplicity $m$ along $Q$.  Let $q=Q\cap C_b$.  Let $\ell-1$ be the order of vanishing of $df_b$ at $q$.

\begin{Cap_Harris_2}
\label{th:Cap_Harris_2}
{\bf \cite[2.6]{CH}}  For any $v\in T_bB$, the image of $\kappa_b(v)$ in $H^0(C_b,\sN_b/(\sN_b)_\tors)$ vanishes to at least order $m-\ell$ at $q$ and cannot vanish to any order $k$ with $m-\ell<k<m$.  If $f(Q)$ is a point, then it vanishes to order at least $m$ at $q$.
\end{Cap_Harris_2}

We combine these results with the existence of a perfect obstruction theory \cite[\S 3.1]{GWinvs} to prove the following.

\begin{gen_smooth}
\label{th:gen_smooth}
Let $e$ be the expected dimension of $\eK_{0,n}(\pd,d,\vec{\varrho})$, which is given by equation \ref{eq:edim}, and assume $e>0$, $d>0$, and $\varrho_i>0$ for all $i$.  The substack $\eK^*_{0,n}(\pd,d,\vec{\varrho})$ is generically smooth of dimension $e$.  At a general point of any irreducible component of $\eK^*_{0,n}(\pd,d,\vec{\varrho})$, the corresponding map $f:C\to\bP^2$ of coarse moduli spaces satisfies the following.
\begin{enumerate}
 \item The multiplicity of $f^*D$ at the $i$th marked point is $\varrho_i$.
 \item There are precisely $(d\delta-\sum\varrho_i)/r$ points where $f^*D$ has multiplicity $r$.
 \item If $e\ge 3$, then $f(C)$ has at worst nodal singularities, none of which lie on $D$.
\end{enumerate}
Moreover, for any irreducible component $\eK\subset\eK^*_{0,n}(\pd,d,\vec{\varrho})$, every evaluation map $\eK\to D$ corresponding to a twisted marking is surjective.
\end{gen_smooth}

\beginproof  Let $\eK\subset\eK^*_{0,n}(\pd,d,\vec{\varrho})$ be an irreducible component, given the reduced induced structure.  Let $B\subset\eK$ be a representable, smooth, dense open substack.  Such a substack exists because the general stable map in $\eK$ has no automorphisms.  Let $b\in B$ be a general point.  Let $\pi:C\to B$ be the coarse moduli space of the universal twisted curve over $B$, and let $f:C\to\bP^2$ be the morphism induced by the universal morphism into $\pd$.  Let $f_b:C_b\to\bP^2$ be the restriction to the fiber over $b$.  Let $t$ be the number of points in the support of $f_b^*D$ and let their multiplicities be $m_i$.

We have the following.
\begin{equation}
\label{eq:inequalities}
e\le\dim T_bB\le 3d-1-\sum(m_i-1)=3d-1-d\delta+t\le e
\end{equation}
The first inequality follows from the general fact that the dimension of a stack is greater than or equal to its expected dimension under a perfect obstruction theory.  For the second inequality, note that $T_bB$ injects into $H^0(C_b,\sN_b)$, that $C_b$ is rational, and that the degree of $\sN_b$ is $3d-2$.  Then apply Lemmas \ref{th:Cap_Harris_1} and \ref{th:Cap_Harris_2} (note that we use $e>0$ here).  The final inequality follows from equation \ref{eq:edim} since the multiplicities of $f^*D$ are bounded below by either the contact types $\varrho_i$ or by $r$.

So each inequality in equation \ref{eq:inequalities} must be an equality.  From the third inequality, we deduce that $t=n+(d\delta-\sum\varrho_i)/r$ which implies statements 1 and 2.  From the second inequality, we deduce that $(\sN_b)_{\tors}=0$.  Combining these facts with Lemma \ref{th:inf_def}, it follows that $\eK^*_{0,n}(\pd,d,\vec{\varrho})$ is generically smooth of dimension $e$.

We have shown that the first order deformations of $f_b:C_b\to\bP^2$ corresponding to vectors in $T_bB$ are the sections of $\sN_b$ which vanish to order $\varrho_i-1$ at a marked point and to order $r-1$ at an unmarked point in $f_b^{-1}D$ (we call these sections \emph{admissible}) and that any such first order deformation extends to a deformation over a curve.  Surjectivity of the evaluation maps $\eK\to D$ now follows from the fact that a section of $\sN_b$ vanishing to order $\varrho_i-1$ at the $i$th marked point cannot fix its image by Lemma \ref{th:Cap_Harris_2}.

For statement 3, since $\sN_b$ has no torsion it suffices to show that at most two points of $C_b$ map to the same point in $\bP^2$, that their tangent directions are distinct, and that at most one point maps to the same point of $D$.  For the latter, if two points map to the same point of $D$, with orders of contact $m_1$ and $m_2$, take any admissible section of $\sN_b$ vanishing to order $m_1-1$ at the first and $m_2$ at the second (we use $e\ge 2$).  The corresponding deformation separates the two points by Lemma \ref{th:Cap_Harris_2}.  For the former, if three points map to the same point of $\bP^2$, then since $e\ge 3$ there is an admissible section vanishing at two of the points and not the third.  Finally, if two points map with the same tangent direction, take any admissible section which vanishes at one of the points but not the other.\qed

\subsection{Maps from nodal twisted curves}
\label{sec:nodal}

In this subsection, we discuss deformations of twisted stable maps which preserve a node.  This is only needed for section \ref{sec:enumerative}, where it is shown that the Gromov-Witten invariants of $\bP^2_{D,2}$, which were computed in \cite{GWinvs}, are enumerative if $D$ is a smooth cubic.

Recall some facts about stable maps into a scheme.  An arbitrary source curve can be described as the result of identifying pairs of points on a smooth proper curve $C$ such that the resulting curve $C'$ is connected.  For any scheme $X$, a morphism $C'\to X$ is the same as a morphism $C\to X$ which sends each pair of identified points to the same point.  In other words, $C'$ satisfies a universal property with respect to $C$.  For curves $C$ over an arbitrary base scheme $S$ together with two disjoint sections $s_1,s_2:S\to C$ in the smooth locus of $C\to S$, there is a clutching construction which produces a curve $C'$ over $S$ and an $S$-morphism $p:C\to C'$ which is universal for morphisms of schemes $f:C\to X$ such that $f\circ s_1=f\circ s_2$ \cite[\S 3]{Kn}.  Because of this, it is clear how to deform a stable map while preserving a node in the source surve.

Such a construction is required for twisted stable maps.  Given a (possibly disconnected) twisted curve $\fC\to S$, an \'etale gerbe $\Sigma\to S$ and two closed embeddings $s_1,s_2:\Sigma\to\fC$ in the smooth locus of $\fC\to S$, it would produce (at least \'etale locally on $S$) a twisted curve $\fC'$ over $S$, an $S$-morphism $p:\fC\to\fC'$, and a 2-morphism $\alpha:p\circ s_1\implies p\circ s_2$ such that the pair $(p,\alpha)$ satisfies a universal property.  To construct $\fC'$, one can use Olsson's description of twisted curves \cite{Ol_twisted}.  So it remains to prove the universal property in some 2-category (which hopefully contains smooth Deligne-Mumford stacks of finite type over $\bC$).

While this is undoubtedly the ``right'' way to approach deformations of twisted curves preserving a node, we adopt a more economical approach here.  Since our application is to genus $0$ stable maps, we only treat disconnecting nodes, and we use that fact that $\pd$ has an open cover by substacks which are global quotients of a scheme by $\bmu_r$ \cite[2.15]{stacks}.
\smallskip

\noindent{\bf Classification of twisted curves.\ } An $n$-marked twisted curve over $\bC$ is determined up to isomorphism by a projective, connected curve $C$ having at worst nodal singularities, a finite set of distinct smooth points $p_1,\ldots,p_n\in C$ (the markings), and a labelling of each marked point and node by a positive integer.  One could construct the corresponding twisted curve by gluing open substacks which \'etale locally look like
\begin{enumerate}
  \item $[\spec \bC[x]/\bmu_r]$, $t\cdot x=t^{-1}x$, at a marked point and
  \item $[\spec \bC[x,y]/(xy)/\bmu_r]$, $t\cdot x=t^{-1}x$, $t\cdot y=ty$, at a node \cite[\S 2.1]{AV}.
\end{enumerate}
This classification also follows from \cite[1.8]{Ol_twisted}.
\smallskip

\noindent{\bf Normalization of a node.\ } Let $\fC$ be a twisted curve over $k$ with coarse moduli space $C$.  For any node $x\in C$, there is a twisted curve $\tilde{\fC}$ which normalizes $x$.  In general, the normalization of a reduced stack is defined by taking a groupoid presentation $R\rightrightarrows U$ and normalizing both $R$ and $U$ \cite[1.18]{Vi}.  If $R\rightrightarrows U$ is a presentation of $\fC$, then $\tilde{\fC}$ is obtained by normalizing only the preimages of $x$ in $R$ and $U$.
\smallskip

We recall some facts about ordinary prestable curves.  Let $B$ be a variety and let $\pi_i:C_i\to B$, for $i=1,2$, be two prestable curves over $B$ with sections $s_i:B\to C_i$ which don't meet the singular locus of the projections $\pi_i$.  Then there is a prestable curve $\pi:C\to B$ and $B$-morphisms $p_i:C_i\to C$ such that
\begin{enumerate}
 \item $p_1\circ s_1=p_2\circ s_2$,
 \item if $X$ is any scheme and $f_i:C_i\to X$ are any morphisms such that $f_1\circ s_1=f_2\circ s_2$, then there is a unique morphism $f:C\to X$ such that$f\circ p_i=f_i$,
 \item and for every geometric point $b\in B$, $\pi_1^{-1}(b)\sqcup\pi_2^{-1}(b)\to\pi^{-1}(b)$ is the normalization of the node $p_1(s_1(b))$.
\end{enumerate}
Moreover, $C$ is unique up to isomorphism, and we say it is obtained by gluing $C_1$ and $C_2$ along $s_1$ and $s_2$.  This follows from \cite[3.4]{Kn}.  Conversely, if $\pi:C\to B$ is a prestable curve such that the singular locus of $\pi$ contains a connected closed subscheme $\Delta\subset C$ which dominates $B$, then there is a finite surjective morphism $\tilde{B}\to B$, two prestable curves $C_1$ and $C_2$ over $\tilde{B}$ and sections $s_i:B\to C_i$ in the smooth loci of the projections, such that $C\times_B\tilde{B}$ is isomorphic to the curve obtained by gluing $C_1$ and $C_2$ along $s_1$ and $s_2$.  This follows from \cite[3.7]{Kn}.
\smallskip

Using these results, we can classify deformations of a twisted stable map that preserve a node.  First, a word about what is meant by this.  Ultimately, we want to know the dimension of an irreducible component of the space of twisted stable maps whose general source curve has a node, in terms of the dimension of the two irreducible components one obtains by separating this node.  So we do not make things very precise, and in fact one must often take a base change or an open covering to go between the two types of deformations.

Let $\fC$ be a connected balanced twisted curve over $\bC$, let $C$ be its coarse moduli space, and let $x\in C$ be a separating node.  Let $F:\fC\to\pd$ be a representable morphism and let $f:C\to\bP^2$ be the induced morphism.  Let $\fC_1$ and $\fC_2$ be the connected components of the normalization of $\fC$ at $x$, let $C_i$ be the coarse moduli space of $\fC_i$, let $x_i\in C_i$ be the preimage of $x$ (viewed as a marked point), let $F_i:\fC_i\to\pd$ and $f_i:C_i\to\bP^2$ be the induced morphisms, and let $\varrho_i$ be the contact type of $F_i$ at $x_i$.  If $f(x)\not\in D$, then the node is untwisted in $\fC$, and it follows from Knudsen's results that a deformation of $F$ which preserves the node $x$ and the condition that $f(x)\not\in D$ is equivalent to a pair of deformations of $F_1$ and $F_2$ which preserve the condition $f_1(x_1)=f_2(x_2)\not\in D$.

Now suppose that $f(x)\in D$.  Then the contact types must be complementary---$\varrho_1+\varrho_2\equiv 0\;(\mathrm{mod\ } r)$---because of the balanced condition.  We claim that a deformation of $F$ which preserves the node $x$ and the condition $f(x)\in D$ is equivalent to a pair of deformations of $F_1$ and $F_2$ which preserve the condition $f_1(x_1)=f_2(x_2)\in D$.  To see this, let $U\subset\bP^2$ be an open set containing $f(x)$ such that $\sO_{\bP^2}(D)$ is trivial on $U$ and let $P\to U$ be an $r$ sheeted cyclic covering which is totally ramified along $D\cap U$ (this is obtained by taking the $r$th root multisection of the tautological section vanishing on $D$).  Then $\mu_r$ acts on $P$ and the stack quotient $[P/\mu_r]$ is canonically identified with $\pd\times_{\bP^2} U$ (see \cite[2.4,2.15]{stacks}).  Given a deformation of $F$ to a family of maps $$(\pi,\tilde{F}):\tilde{\fC}\to B\times\pd$$ which has a closed substack $\Delta\subset\tilde{\fC}$ contained in the singular locus of $\pi$ and containing $x$, let $V\subset\tilde{\fC}$ be the preimage of $U$.  The fiber product $V\times_{\pd}P$ is a curve over $B$ since $\tilde{F}$ is representable.  One obtains the deformations of $F_1$ and $F_2$ by normalizing the preimage of $\Delta$ in $V\times_{\pd}P$ and taking stack quotients by $\bmu_r$.

Conversely, if we are given deformations of $F_1$ and $F_2$ over a variety $B$, then we take analogous fiber products and simply need to glue the resulting curves to form nodes.  So we have curves $\tilde{C}_i$ over $B$ and subschemes $\sigma_i\subset\tilde{C}_i$.  These subschemes are finite and \'etale over $B$ of degree equal to the greatest common divisor of $r$ and $\varrho_1$.  In fact, they are the total spaces of principal bundles over $B$ having cyclic structure group.  So after an \'etale base change on $B$ we can assume they are trivial.  On the fibers over the special point $b\in B$ (corresponding to the original map $F$) we are given an identification of $\sigma_1$ with $\sigma_2$.  This extends uniquely to an isomorphism $\sigma_1\to\sigma_2$ over $B$.  Thus the curves can be glued to a curve $\tilde{C}$ using Knudsen's construction.  Since we assumed that the condition $f(x_1)\in D$ is preserved by the deformation, and since $P\to U$ is totally ramified over $D$, the morphisms $\tilde{C}_i\to P$ extend to a morphism $\tilde{C}\to P$.  Now taking the stack quotient of $\tilde{C}$ yields the required deformation of $F$.

\section{Intersection numbers on \boldmath $\eK_{0,n}(\pd,d,\vec{\varrho})$}
\label{sec:cap_harris}
\setcounter{equation}{0}

In this section we prove an analogue of a formula of Caporaso and Harris for rational plane curves having prescribed tangencies with a smooth plane curve of positive genus.  Throughout this section we fix a plane curve $D$ of degree $\delta$, positive integers $r$ and $d$, and an $n$-tuple of integers $\varrho_1,\ldots,\varrho_n$ such that $1\le\varrho_i\le n$ and $\sum\varrho_i = d\delta$.  Let $\eK=\eK_{0,n}(\pd,d,\vec{\varrho})$.

Let $A^*(\eK)$ to denote the operational Chow ring of $\eK$ \cite{Vi}.  We define $N^1(\eK)$ to be $A^1(\eK)$ modulo the equivalence relation $a_1\equiv a_2$ if for all $b\in A_1(\eK)$, $$\int_{\eK}a_1\cap b=\int_{\eK}a_2\cap b.$$  Let $\eC$ be the universal coarse curve over $\eK$.  By this we mean that we have a representable morphism $\eC\to\eK$ such that for every morphism $B\to\eK$ from a scheme $B$, the pullback of $\eC$ to $B$ is canonically isomorphic to the coarse moduli space of the pullback of the universal twisted curve to $B$.  So we have the following diagram.
$$\xymatrix{
\eC \ar[r]^f \ar[d]^{\pi} & \bP^2 \\
\eK}$$

We define classes $h,\chi_1,\ldots,\chi_n\in N^1(\eK)$ as follows.  Let $h=\pi_*(f^*p)$, where $p\in A^2(\bP^2)$ is the class of a point and $\pi_*$ is the proper, flat pushforward.  In the notation of \cite[\S 17]{Fu}, this equals $\pi_*(f^*p\cdot [\pi])$, where $[\pi]$ is the orientation class.  Let $\chi_i=e_i^*\tilde{p}$, where $e_i:\eK\to D$ is the $i$th evaluation map and $\tilde{p}\in N^1(D)$ is the class of a point.

The following relation between these classes leads to the analogue of Caporaso and Harris's formula.

\begin{relations}
\label{th:relations}
If $r>d\delta$, then $$h=\sum_{i=1}^n \varrho_i\chi_i.$$
\end{relations}

\beginproof  First we derive the key consequence of the hypothesis $r>d\delta$ and in doing so fix some notation.  Let $F:\fC\to\pd$ be a twisted stable map over $\bC$ which lies in $\eK$, and let $C$ be the coarse curve of $\fC$, and let $g:C\to\bP^2$ be the induced morphism.  Let $x_1,\ldots,x_n\in\fC$ be the marked points and let $y_1,\ldots,y_m$ be the twisted nodes of $\fC$ which lie on at least one component of $\fC$ which maps with positive degree.  By renumbering the markings if necessary, assume that $x_1,\ldots,x_k$ lie on components of $\fC$ which map with positive degree and assume that $x_{k+1},\ldots,x_n$ lie on components which map with degree $0$.

First note that no twisted node $y_i$ can lie on two components which map with positive degree.  Indeed, since the contact types of $y_i$ on the two components must add to $r$, it would follow that the pullback of $D$ to the normalization of $C$ has degree at least $r$, contradicting $r>d\delta$.  So each $y_i$ lies on a unique component mapping with positive degree.  Let $\sigma_i$ be the contact type of $y_i$ on this component.  For the same reason, the contact types $\sigma_i$ for all $i$ and $\varrho_j$ for $1\le j\le k$ are equal to the intersection number between $D$ and the component of $C$ at the given point.

Suppose that $A\subset\{k+1,\ldots,n\}$ is the set of markings lying in a fixed connected component of $g^{-1}(D)$ and suppose that $B\subset\{1,\ldots,m\}$ is the set of nodes lying in the same component.  We claim that 
\begin{equation}
\label{eq:multiplicities}
\sum_{i\in B}\sigma_i=\sum_{j\in A}\varrho_j.
\end{equation}
That they are congruent modulo $r$ follows from the fact that the sum of contact types at all twisted points lying on an irreducible component mapping with degree $0$ must be a multiple of $r$, together with the fact that the two contact types at a node sum to a multiple of $r$.  From this congruence, equality is deduced from the fact that both sides are between $0$ and $r-1$.

To show that $h=\sum\varrho_i\chi_i$, it suffices to show that for any one dimensional integral closed substack $V\subset\eK$, $\int_V h=\sum\varrho_i\int_V\chi_i$.  It also suffices to replace $V$ with its normalization.  Let $\eC_V$ be the restriction of the universal curve to $V$.  Then $\int_V h=\int_{\eC_V} f_V^* p = \int_{\bP^2} p\cap (f_V)_*[\eC_V]$ which is the degree of $f_V$.  There exists a finite, flat base change $W\to V$ such that every irreducible component of $\eC_W$ is generically irreducible over $W$ (hence generically smooth).  Replace $V$ with $W$.

Now we use the notation introduced at the beginning of the proof, where we take $F:\fC\to\pd$ to be the map corresponding to a general point of $V$.  For each node $y_i$, there are two irreducible components of $\eC_V$ which contain $y_i$, and the intersection between these components is a section $t_i:V\to\eC_V$.  Let $s_i:V\to\eC_V$ be the section corresponding to the $i$th marking.  The degree of $f_V$ can be computed by $(f_V)_*(f_V)^*[D]=(\deg(f_V))[D]$.  This shows that $$\deg(f_V)=\sum_{i=1}^k\varrho_i\deg(f_V\circ s_i)+\sum_{i=1}^m\sigma_i\deg(f_V\circ t_i).$$  Here the morphisms on the right hand side are viewed as going from $V$ to $D$.  The first summation is $\sum_{i=1}^k\varrho_i\int_V\chi_i$.  It follows from equation \ref{eq:multiplicities} that the second summation is $\sum_{i=k+1}^n\varrho_i\int_V\chi_i$.\qed

Now we define intersection numbers on the stack $\eK^*_{0,n}(\pd,d,\vec{\varrho})$.  First we introduce some notation.  If $\alpha=(\alpha_1,\ldots)$ is a sequence of positive integers with all but finitely many equal to $0$, then let $\vert\alpha\vert=\sum\alpha_i$, $I\alpha=\sum i\alpha_i$, and $\alpha!=\prod\alpha_i!$.

\begin{CH_numbers}
\label{th:CH_numbers}
Let $\alpha=(\alpha_1,\ldots)$ and $\beta=(\beta_1,\ldots)$ be sequences of positive integers with all but finitely many equal to $0$.  Let $n=\vert\alpha\vert + \vert\beta\vert$, let $\vec{\varrho}$ be an $n$-tuple of positive integers such that $\varrho_j=i$ for $\alpha_1+\ldots+\alpha_{i-1}<j\le\alpha_1+\ldots+\alpha_i$ and for $I\alpha+\beta_1+\ldots+\beta_{i-1}<j\le I\alpha+\beta_1+\ldots+\beta_i$, and choose $r$ so that $r>d\delta$.  Let $D$ be a smooth curve of degree $\delta$ such that $d\delta=I\alpha+I\beta$.  If the expected dimension of $\eK^*_{0,n}(\bP^2_{D,r},d,\vec{\varrho})$ equals $e\ge\max(1,\vert\alpha\vert)$, then we define $$N_d^D(\alpha,\beta)=\frac{1}{\beta!}\int_{\eK^*_{0,n}(\pd,d,\vec{\varrho})} h^{e-\vert\alpha\vert}\cdot\prod_{i=1}^{\vert\alpha\vert}\chi_i.$$
\end{CH_numbers}

These numbers all have enumerative significance according to the following proposition (hence are independent of $r$).  When we say that a curve has a specified $i$th order contact with $D$, we mean that we have chosen a general point of $D$ and require the curve to have an $i$th order contact with $D$ at this point.

\begin{CH_enum}
\label{th:CH_enum}
$N_d^D(\alpha,\beta)$ is the number of rational degree $d$ curves passing through $e-I\alpha$ general points, having $\alpha_i$ specified $i$th order contacts with $D$, and having $\beta_i$ unspecified $i$th order contacts with $D$.  If $e\ge 3$ or $d\le 2$, then these curves have at worst nodal singularities, none of which lie on $D$.
\end{CH_enum}

\beginproof Let $\eC^{e-I\alpha}$ be the $(e-I\alpha)$-fold product of the universal curve with itself.  We have a morphism $$F:\eC^{e-I\alpha}\to (\bP^2)^{e-I\alpha}\times D^{\vert\alpha\vert}$$ given by evaluation at marked points.  The integral in Definition \ref{th:CH_numbers} is equal to the degree of $F$.  We can throw away the boundary of $\eC^{e-I\alpha}$ to get a proper morphism of smooth schemes $U\to V$, whose degree equals the integral in question.  Given a curve satisfying the stated conditions, there are $\beta!$ such twisted stable maps, because the marked points corresponding to unspecified contacts can be relabeled arbitrarily.  Now the proposition follows from generic smoothness, together with Theorem \ref{th:gen_smooth}.\qed

Now we come to the main theorem of this section.  Let $e_k$ be the sequence with all but the $k$th entry $0$, and the $k$th entry equal to $1$.

\begin{CH_recursion}
\label{th:CH_recursion}
If $e-I\alpha>0$, then $$N_d^D(\alpha,\beta)=\sum_{k:\beta_k>0} kN_d^D(\alpha+e_k,\beta-e_k).$$
\end{CH_recursion}

\beginproof If $j$ is chosen so that $\varrho_j=k$ and $j>I\alpha$, then $$N_d^D(\alpha+e_k,\beta-e_k)=\frac{\beta_k}{\beta!}\int_{\eK^*_{0,n}(\pd,d,\vec{\varrho})}h^{e-I\alpha-1}\cdot\prod_{i=1}^{I\alpha}\chi_i\cdot\chi_j.$$  Therefore, it suffices to show that in $N^1(\eK^*_{0,n}(\pd,d,\vec{\varrho}))$, we have $h=\sum_{i=1}^n \varrho_i\chi_i$ and $\chi_i^2=0$.  The latter equality is obvious, because if one chooses two distinct points of $D$, then their preimages under $e_i:\eK^*_{0,n}(\pd,d,\vec{\varrho})\to D$ are disjoint.  The former equality follows from Proposition \ref{th:relations}.\qed

\section{Derivation of the main result}
\label{sec:kontsevich}
\setcounter{equation}{0}

Let $E\subset\bP^2$ be a smooth cubic.  Let $d$ be a positive integer and let $a,b,c$ be nonnegative integers such that $a+2b+2c\le 3d$ and $a+2b+c\le 3d-1$.  We define $N_d(a,b,c)$ to be $N_d^E((a,b),(3d-a-2b-2c,c))$.  As a convention, we'll say that $N_d(a,b,c)=0$ if any of the above hypotheses are not satisfied.  Our recursion for these numbers will show that they are independent of $E$.

From Theorem \ref{th:CH_recursion} we have
\begin{equation}
\label{eq:CH}
N_d(a,b,c)=N_d(a+1,b,c)+2N_d(a,b+1,c-1) \mbox{ if } a+2b+c<3d-1 \mbox{ and } a,b\ge 0.
\end{equation}
To apply this formula, we define a family of functions $P_d^{b,c}:\bZ\to\bZ$ by $$P_d^{b,c}(t)=N_d(3d-1-2b-c-t,b,c).$$  The parameter $t$ is equal to the exponent of $h$ in Definition \ref{th:CH_numbers}, so it equals the number of general points in $\bP^2$ through which the curves are required to pass.  This function vanishes outside of the interval $\max(0,c-1)\le t\le 3d-1-2b-c$.

Let $\Delta$ be the difference operator defined by $\Delta P(t) = P(t+1)-P(t)$.  Then equation \ref{eq:CH} can be rewritten as
\begin{equation}
\label{eq:CH2}
\Delta P_d^{b,c}(t) = 2P_d^{b+1,c-1}(t) \mbox{ if } 0\le t \le 3d-2-2b-c \mbox{ and } b\ge 0.
\end{equation}

This equation implies that the restriction of $P_d^{b,c}$ to the interval $0\le t\le 3d-1-2b-c$ agrees with an integer polynomial $Q_d^{b,c}$ of degree at most $c$.  If $2b+2c\neq 3d$, then there are at least $c+1$ integers in this interval, so $Q_d^{b,c}$ is uniquely determined by $P_d^{b,c}$.  If $2b+2c=3d$, then equation \ref{eq:CH} becomes $N_d(0,b,c)=2N_d(0,b+1,c-1)$, so $N_d(0,b,c)=2^{c-1}N_d(0,b+c-1,1)$, which translates into $$Q_d^{b,c}(c-1)=2^{c-1}Q_d^{b+c-1,1}(0).$$  The general solution for $Q_d^{b,c}(t)$ below satisfies this equality.

Since $Q_d^{b,c}(t)$ has degree at most $c$, we can write
$$Q_d^{b,c}(t)=\alpha_0\ch{t}{c} + \alpha_1\ch{t}{c-1} +\cdots + \alpha_c.$$  Since $Q_d^{b,c}(t)=0$ for $0\le t\le c-2$, we have $\alpha_i=0$ for $i\ge 2$.  We claim that $\alpha_1=2\alpha_0$.  Since $\Delta\ch{t}{c}=\Delta\ch{t}{c-1}$, equation \ref{eq:CH2} reduces to the case $c=1$.  In that case $Q_d^{b,1}(t)=\alpha_0t+\alpha_1$ and $Q_d^{b+1,0}(t)=\alpha_0/2$.  So we have reduced to showing that $Q_d^{b,1}(0)=4Q_d^{b+1,0}(0)$.  Translating into the notation of Definition \ref{th:CH_numbers} and replacing $b$ with $b+1$, this says the following.

\begin{key_relation}
\label{th:key_relation}
If $a+2b=3d$, then $$N_d^E((a,b-1),(0,1))=4N_d^E((a-1,b),(1,0)).$$
\end{key_relation}

\beginproof  Let $E^{a+b}\to\Pic^{3d}E$ be the morphism $$(x_1,\ldots,x_a,y_1,\ldots,y_b) \mapsto \sO_E(\sum x_i+2\sum y_i).$$  Let $S\subset E^{a+b}$ be the preimage of $[\sO_E(d)]$, where $\sO_E(1)$ is the restriction of $\sO_{\bP^2}(1)$ to $E$.  Let $p_i:S\to E$, $1\le i\le a$ and $q_i:S\to E$, $1\le i\le b$ be the projections onto the coordinates $x_i$ and $y_i$ respectively.  Let $\hat{p}_i:S\to E^{a+b-1}$ and $\hat{q}_i:S\to E^{a+b-1}$ be the projections onto the complementary factors.  Then the following diagrams are fiber squares.
$$\xymatrix{
S \ar[r]^{\hat{p_a}} \ar[d]_{p_a} \ar@{}[dr]|{\Box} & E^{a+b-1} \ar[d] \\
E \ar[r]^{\cong} & \Pic^1(E)}$$
Here the right arrow is the composition of $E^{a+b-1}\to\Pic^{3d-1}(E)$ with the morphism $\Pic^{3d-1}E\to\Pic^1(E)$ sending $[\sL]$ to $[\sO_E(d)\otimes\sL^{-1}]$.
$$\xymatrix{
S \ar[r]^{\hat{p}_b} \ar[d]_{p_b} \ar@{}[dr]|{\Box} & E^{a+b-1} \ar[d] \\
E \ar[r] & \Pic^2(E)}$$
Here the right arrow is defined similarly to above and the bottom arrow is the squaring map.  Letting $S\to E_{a+b}$ be the projection, it is now clear that $N_d^E((a-1,b),(1,0))$ equals the number of rational degree $d$ curves passing through a general divisor in the image of $S$.  Since we have shown that $\hat{p}_b$ is an \'etale degree 4 cover, it follows that $N_d^E((a,b-1),(0,1))$ is the number passing through any of 4 such divisors (while the divisors are not general with respect to each other, they can all be chosen to lie in a given dense open subset of $S$).\qed

From all this, it follows that there exist numbers $K^{\lambda}_d$ such that
\begin{equation}
\label{eq:main_1}
N_d(3d-1-2b-c-t,b,c)=2^cK^{b+c}_d\left[\ch{t}{c}+2\ch{t}{c-1}\right].
\end{equation}
In fact, $K^{b+c}_d$ is an integer whenever $2b+2c\neq 3d$, since then $K^{b+c}_d=N_d(0,b+c,0)$.  If $2b+2c=3d$, then $K^{b+c}_d$ must be in $\frac{1}{4}\bZ$.

In the next section we show that the numbers $N_d(a,0,c)$ equal the Gromov-Witten invariants of $\bP^2_{E,2}$ which were computed in \cite{GWinvs}.  From that computation we obtain the base cases $$K_1^0=K_1^1=1,\; K_2^3=\frac{3}{4},$$ and recursion 5.7 from [ibid] leads to the following recursion for the coefficients $K_d^b$.
\footnotesize
\begin{eqnarray}
\lefteqn{K_d^{(b)}f^{(b)}_d = \sum_{\stackrel{\scriptstyle b_1+b_2=b}{d_1+d_2=d}}K_{d_1}^{(b_1)}K_{d_2}^{(b_2)}f^{(b_1)}_{d_1}f^{(b_2)}_{d_2}\left[d_1^2d_2^2\ch{3d-4-b}{3d_1-2-b_1}-d_1^3d_2\ch{3d-4-b}{3d_1-1-b_1}\right] +} \label{eq:main_2} \\
& & \sum_{\stackrel{\scriptstyle b_1+b_2=b-1}{d_1+d_2=d}}K_{d_1}^{(b_1)}K_{d_2}^{(b_2)}f^{(b_1)}_{d_1}f^{(b_2)}_{d_2}\alpha^{\vec{b}}_{\vec{d}}\left[2d_1d_2\ch{3d-4-b}{3d_1-2-b_1}-d_1^2\ch{3d-4-b}{3d_1-1-b_1}-d_2^2\ch{3d-4-b}{3d_1-3-b_1}\right], \nonumber
\end{eqnarray}
\normalsize
where 
\begin{equation}
\label{eq:main_3}
\alpha^{\vec{b}}_{\vec{d}}=\frac{(3d_1-2b_1)(3d_2-2b_2)(3d_2-1)}{3d_2(3d_2-1-b_2)}.
\end{equation}
In the sums, $d_i$ are positive integers and $b_i$ are nonnegative integers.  The recursion is only valid for values of $b$ and $d$ not covered by the base cases, and we set $f_d^{(b)}$ and $\alpha_{\vec{d}}^{\vec{b}}$ to be 0 whenever their denominators would make them undefined.

\section{Enumerativity of certain Gromov-Witten invariants}
\label{sec:enumerative}
\setcounter{equation}{0}

Let $p\in N^2(\bP^2)$ be the class of a point, $\alpha\in N^0(E)$ be the fundamental class, and $\beta\in N^1(E)$ be the class of a point.  We use the notation $I_d(p^k\alpha^{\ell}\beta^m)$ to denote genus $0$ Gromov-Witten invariants of $\bP^2_{E,2}$, as in \cite{GWinvs}.  This section is devoted to proving the following.

\begin{enumerative}
\label{th:enumerative}
Given integers $a,c,d$ such that $d>0$, $a,c\ge 0$, $a+2c\le 3d$, and $a+c\le 3d-1$, $$N_d(a,0,c)=\frac{1}{(3d-a-2c)!}I_d(p^{3d-1-a-c}\alpha^{3d-a-2c}\beta^a).$$  Moreover, this is the number of rational degree $d$ plane curves passing through $3d-1-a-c$ general points of $\bP^2$, passing through $a$ general points of $E$, and having $c$ tangencies to $E$.  These curves have at worst nodal singularities, none of which lie on $E$.
\end{enumerative}

It is convenient to replace the contact type $\vec{\varrho}$ with the number of twisted marked points, and assume the markings are ordered so that the twisted markings come before the untwisted ones.  So let $\eK_{0,n}(\pe,d,m)=\eK_{0,n}(\pe,d,(1,\ldots,1,0,\ldots,0))$, where there are $m$ ones and $n-m$ zeros.  The expected dimension of $\eK_{0,n}(\pe,d,n)$ is $$n+\frac{3d-n}{2}-1.$$ The number of tangencies being imposed is $(3d-n)/2$, and must be an integer.  Therefore, the expected dimension is at least equal to $1$ and is at least $3$ if $d\ge 3$.

Recall that Gromov-Witten invariants are defined to be integrals over $\eK_{0,n}(\pe,d,m)$ \cite[\S 2.3]{GWinvs}.  By Theorem \ref{th:gen_smooth}, the closed substack $\eK^*_{0,n}(\pe,d,m)$ has the expected dimension, and therefore the restriction of the virtual fundamental class to $\eK^*_{0,n}(\pe,d,m)$ is the ordinary fundamental class.  It follows from the definitions that the contribution of $\eK^*_{0,n}(\pe,d,m)$ to the Gromov-Witten invariants gives precisely $N_d(a,0,c)$.  In light of Proposition \ref{th:CH_enum}, it remains only to show that the remaining irreducible components of $\eK_{0,n}(\pe,d,m)$ do not contribute to the invariants.  For this we need the following lemmas.

\begin{smooth_source}
\label{th:smooth_source}
If $\eK\subset\eK_{0,n}(\pe,d,n)$ is an irreducible component whose general point corresponds to a stable map $\ff:\fC\to\pe$ with $\fC$ smooth, then the associated map of coarse moduli spaces $f:C\to\bP^2$ maps $C$ birationally onto its image.
\end{smooth_source}

\beginproof  Suppose that $C$ does not map birationally onto its image.  Let $\bar{C}$ be the normalization of the image of $f$, so that $f$ decomposes into
$$\xymatrix{
C \ar[r]^g & \bar{C} \ar[r]^h & \bP^2.}$$
The morphism $h$ corresponds to a representable morphism $\bar{\fC}\to\pe$ from a smooth twisted curve $\bar{\fC}$ which on the level of coarse moduli spaces equals $h$.  The point of $\fC$ lying over a point $x\in C$ is twisted if and only if $k:=\mult_x f^*E$ is odd.  Likewise, the point of $\bar{C}$ lying over $g(x)$ is twisted if and only if $\ell:=\mult_{g(x)}h^*E$ is odd.  Since $g$ is ramified to order $k/\ell$ at $x$, it follows that $g$ lifts uniquely to a representable morphism $\fC\to\bar{\fC}$ which factors $f$.

Let $e$ be the degree of $h$ and let $m$ be the number of twisted marked points on $\bar{\fC}$.  We have shown that $h:\bar{\fC}\to\pe$ lies on an irreducible component of $\eK_{0,m}(\pe,e,m)$ which has the expected dimension $(3e+m)/2-1$.  Moreover, the dimension of $\eK$ is at least equal to the expected dimension, which is $(3d+n)/2-1$.  Since we assumed that $f:\fC\to\pe$ was chosen generally, we must have $3e+m\ge 3d+n$.  Since $m\le 3e$, we have $6e\ge 3d$, so $g$ must be a double cover and $6e=3d$, so $m=3e$ and $n=0$.  Since $C$ and $\bar{C}$ are rational, $g$ is ramified at exactly 2 points.  But $g$ must be ramified over every twisted point of $\bar{\fC}$, which is impossible.\qed

For each untwisted marking we have an evaluation map $\eK_{0,n}(\pe,d,m)\to\bP^2$, and for each twisted marking an evaluation map $\eK_{0,n}(\pe,d,m)\to E$.  In the following proof, we use the result of section \ref{sec:nodal}, which compares deformations of a twisted stable map to deformations of the pair of maps obtained by normalizing a separating node.

\begin{gen_elmt}
\label{th:gen_elmt}
Let $\eK\subset\eK_{0,n}(\pe,d,m)$ be an irreducible component.  Let $f:\fC\to\pe$ be a general map in $\eK$ and let $\tau$ be the number of nodes of $\fC$.  Let $ev:\eK\to (\bP^2)^{n-m}\times E^m$ be the product of all evaluation maps.  Then $$\dim(ev(\eK))+\tau\le\edim(\eK).$$  Moreover, each evaluation map is surjective, except possibly for evaluation maps corresponding to untwisted markings which lie on a component of $\fC$ mapping into $E$ with degree $0$.  In the latter case, the evaluation map surjects onto $E$.
\end{gen_elmt}

\beginproof  We prove this by induction on $\tau$.  First assume $\tau=0$.  If $d=0$, then a simple calculation shows that $\edim(\eK)\ge 1$ and $\edim(\eK)\ge 2$ if there are no twisted markings.  The surjectivity of evaluation maps in this case is clear.  If $d>0$, then Lemma \ref{th:smooth_source} implies that $\eK\subset\eK^*_{0,n}(\pe,d,m)$.  So the result follows from Theorem \ref{th:gen_smooth} if it is shown that untwisted evaluation maps $\eK\to\bP^2$ surject.  If not, then the image is at most 1-dimensional, so by forgetting all untwisted markings we would have a component of $\eK^*_{0,m}(\pe,d,m)$ of dimension $0$.  But this contradicts the fact that the expected dimension is $(3d+m)/2-1\ge 1$.

For the general case, note first that $\fC$ has no nodes between components which map with degree $0$.  If such a node maps to a point not in $E$, then this follows from the irreducibility of $\bar{M}_{0,n}$.  If the node maps into $E$, then it follows from the fact that $\eK_{0,n}(B\mu_2)$ is flat over $\bar{M}_{0,n}$ \cite[3.0.5]{ACV}.

Suppose that $\fC$ has a node $x$ lying on two components which map with positive degree.  Let $\fC_1$ and $\fC_2$ be the twisted curves resulting from normalizing $x$, and let $\eK_1$ and $\eK_2$ be the irreducible components of the stack of twisted stable maps in which $f\vert_{\fC_1}$ and $f\vert_{\fC_2}$ lie (we take the preimages of $x$ to be marked points).  By induction, the result holds for $\eK_1$ and $\eK_2$.  If $x$ is untwisted, then $ev(\eK)=ev(\eK_1)\times_{\bP^2}ev(\eK_2)$, where the morphisms to $\bP^2$ are projections onto the factor corresponding to the preimage of $x$.  Since these evaluation maps are surjective, it follows that $\dim(ev(\eK))=\dim(ev(\eK_1))+\dim(ev(\eK_2))-2$.  Moreover, $\edim(\eK)=\edim(\eK_1)+\edim(\eK_2)-1$, so it follows that $\dim(ev(\eK))+\tau\le\edim(\eK)$.  The surjectivity of evaluation maps also follows by induction.

If $x$ is twisted, then replacing $\bP^2$ with $E$ we find that $\dim(ev(\eK))=\dim(ev(\eK_1))+\dim(ev(\eK_2))-1$, and a calculation shows that $\edim(\eK)=\edim(\eK_1)+\edim(\eK_2)$.  So the result follows in the same way.

The other case to consider is when the node $x$ joins a positive degree component with a degree $0$ component.  Then the above argument will fail if the degree $0$ component maps into $E$ and the node is untwisted.  However, in this case the condition that the preimage of $x$ in $\fC_1$ maps to $E$ imposes a nontrivial condition on the space $\eK_1$, and since the evaluation maps at least surject onto $E$, it follows that $\dim(ev(\eK))=\dim(ev(\eK_1))+\dim(ev(\eK_2))-2$ anyway.  Everything else works as before.\qed

We now return to the proof of Theorem \ref{th:enumerative}.  If $\eK\subset\eK_{0,n}(\pe,d,m)$ is any irreducible component which is not in $\eK^*_{0,n}(\pe,d,m)$, then by Lemma \ref{th:smooth_source} we have $\tau>0$ in the notation of Lemma \ref{th:gen_elmt}.  Therefore, $\dim(ev(\eK))<\edim(\eK)$ for all irreducible components not in $\eK^*_{0,n}(\pe,d,m)$.  But this implies that the pushforward of the virtual fundamental class of $\eK_{0,n}(\pe,d,m)$ under $ev$ equals the pushforward of its restriction to $\eK^*_{0,n}(\pd,d,m)$ (which is well-defined, since the latter has the expected dimension).  This implies that irreducible components not in $\eK^*_{0,n}(\pe,d,m)$ do not contribute to the Gromov-Witten invariants.

\noindent{\bf Remark}  There exist irreducible components of $\eK_{0,n}(\pe,d,m)$ having greater than the expected dimension.  One way to see this is to consider maps from a twisted curve with two irreducible components, one of which maps with degree $0$ and contains at least $5$ twisted points.

\noindent University of Michigan \\
2074 East Hall \\
Ann Arbor, MI 48109-1043 \\
\ttfamily{cdcadman@umich.edu}
\end{document}